\documentclass[a4paper,10pt]{article}
\usepackage[utf8]{inputenc}
\usepackage{amsmath,amsthm}
\usepackage{amsfonts}
\usepackage{amssymb}
\usepackage{a4wide}
\usepackage{todonotes}
\usepackage{tikz}
\usepackage{xcolor}
\usetikzlibrary{patterns,patterns.meta}

\newcommand{\R}{\mathbb{R}}
\newcommand{\N}{\mathbb{N}}
\newcommand{\Z}{\mathbb{Z}}

\newcommand{\norm}[2]{\left\Vert #1 \right\Vert_{#2}}

\DeclareMathOperator*{\essinf}{essinf}
\DeclareMathOperator{\ran}{ran}

\newcommand{\X}{\mathcal{X}}
\newcommand{\Y}{\mathcal{Y}}

\renewcommand{\d}{\,\mathrm d}
\newcommand{\leb}{\mu_{\mathrm{Leb}}}

\newtheorem{claim}{Claim}
\newtheorem{lem}{Lemma}
\newtheorem{theorem}{Theorem}
\newtheorem{ex}{Example}
\newtheorem{Def}{Definition}
\newtheorem{ass}{Assumption}
\newtheorem{pro}{Proposition}
\newtheorem{rem}{Remark}
\newtheorem{cor}{Corollary}

\title{A unified concept of the degree of ill-posedness for compact and non-compact linear operator {equations} in Hilbert spaces under the auspices of the spectral theorem}
\author{Frank Werner \footnote{Institute for Mathematics, University of W\"urzburg, Emil-Fischer-Str.~30, 97074~W\"urzburg, Germany. Email: frank.werner@uni-wuerzburg.de}
\and Bernd Hofmann \footnote{Chemnitz University of Technology, Faculty of Mathematics, 09107 Chemnitz, Germany. Email: hofmannb@mathematik.tu-chemnitz.de}}

\begin{document}

\maketitle

\begin{abstract}
Covering ill-posed problems with compact and non-compact operators regarding the degree of ill-posedness is a never ending story written by many authors in the inverse problems literature. This paper tries to add a new narrative and some new facets with respect to this story under the auspices of the spectral theorem. The latter states that any self-adjoint and bounded operator is unitarily equivalent to a multiplication operator on some (semi-finite) measure space. We will exploit this fact and derive a distribution function from the corresponding multiplier, the growth behavior of which at zero allows us to characterize the degree of ill-posedness. We prove that this new concept coincides with the well-known one for {equations with} compact operators (by means of their singular values), and illustrate the implications along examples with non-compact operators, including the Hausdorff moment operator and convolutions.
\end{abstract}

\section{Introduction}

Let $T : \X \to \Y$ be a bounded and injective linear operator with range ${\rm ran}(T)$, which is mapping between infinite dimensional separable real Hilbert spaces $\X$ and $\Y$. In this study, we consider the operator equation
\begin{equation}\label{eq:model}
Tu = g
\end{equation}
that is ill-posed in correspondence with the fact that ${\ran}(T)$ is assumed to be a non-closed subset of $\Y$. In the seminal paper \cite{Nashed86}, the case of {\sl compact operators} $T$ well-studied in the literature
 was called {\sl ill-posedness of type~II}, and their strength of ill-posedness is clearly expressed by the decay rate of the singular values of $T$. The majority of occurring linear ill-posed problems
 (e.g., Fredhom and Volterra integral equations of first kind over some bounded domain in $\R^d$ as well as history match problems for the heat equation with homogeneous boundary conditions) belong to that type.
 Such problems with compact forward operators are extensively analyzed in the inverse problems literature, see e.g. \cite{h00,LS06,n86}. In the Hilbert space setting, the alternative {\sl ill-posedness of type~I} is characterized by
{\sl non-compact operators} $T$ with non-closed range and arises for convolution operators on unbounded domains (see Subsection~\ref{sec:conop} below) as well as for Ces\`{a}ro operators (cf., e.g.,~\cite{Brown65,DFH24,KinHof24}) and in the context of the Hausdorff moment problem (see Subsection~\ref{sec:hausop} below). The absence of singular values makes the things more difficult for non-compact $T$. {Therefore, it is our goal to make a suggestion for classifying the strength of ill-posedness also for non-compact operators $T$ with non-closed range by means of their spectral properties.}

{Covering ill-posed problems with compact and non-compact operators regarding the nature and degree of ill-posedness is a never ending story written by many authors in the inverse problems literature. For example, we refer to the claim of Nashed in \cite{Nashed86} that {equations with} non-degenerating compact operators mapping between infinite dimensional Hilbert spaces are always {\sl more ill-posed} than {those with} non-compact ones. Arguments for that claim come from the comparison of ill-posed situations with two linear operators by partial orderings and associated range inclusions (see, e.g, \cite{HofKind10,MH24}).
In this context, the range of a non-compact operator can never be a subset of the range of a compact operator. Also in \cite{HofKind10}, however, it was shown by observing the modulus of injectivity for various discretizations that the ill-posedness nature of non-compact operators is completely different from the nature of compact ones and that hence the simple claim of Nashed is rather problematic, because sequences of compact operators do never converge in norm to a non-compact operator.
In the present paper, we try to discuss some \emph{new facets} of the story of covering compact and non-compact cases under the auspices of the spectral theorem.}

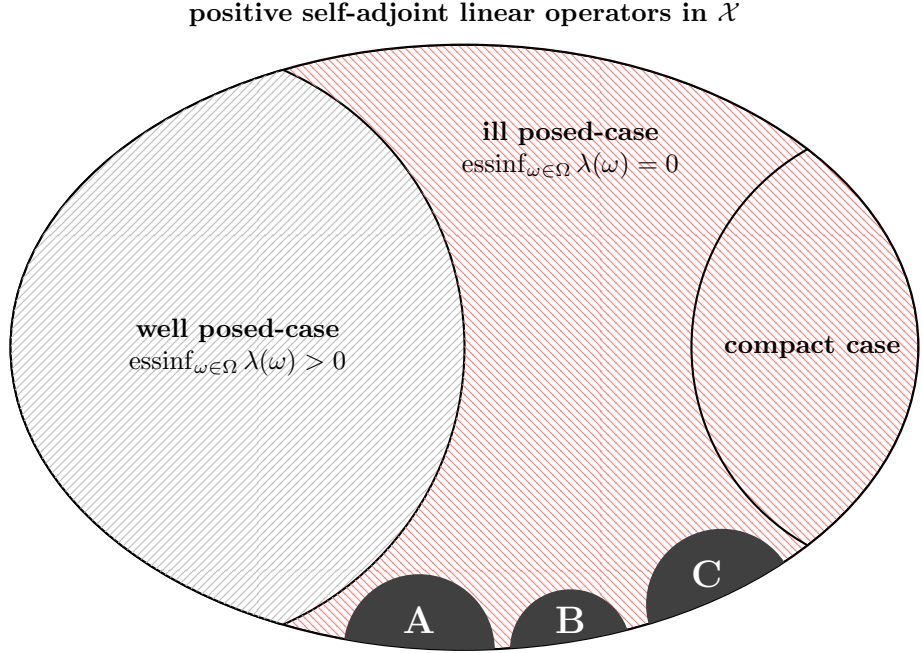
\begin{figure}[!htb]
	
\begin{tikzpicture}[scale = 2,thick]
	
	\draw[black] (0,0) ellipse (3cm and 2cm);
	\node[align=center] at (0,2.2) {\textbf{positive self-adjoint linear operators in $\X$}};
	\begin{scope}
	\clip (0,0) ellipse (3cm and 2cm);
	\draw[pattern=north east lines,pattern color = lightgray] (-2,0) circle(2cm);
	\end{scope}
	\node[align=center] at (-1.5,0) {\textbf{well posed-case} \\$\essinf_{\omega \in \Omega} \lambda(\omega) > 0$};
	\begin{scope}[even odd rule]
		\clip (0,0) ellipse (3cm and 2cm) 
		(-2,0) circle(2cm); 
		\draw[pattern=north west lines,pattern color = red!50!white] (0,0) ellipse (3cm and 2cm);		
	\end{scope}
	\node[align=center] at (.7,1.3) {\textbf{ill posed-case}\\$\essinf_{\omega \in \Omega} \lambda(\omega) = 0$};
	\begin{scope}
		\clip (0,0) ellipse (3cm and 2cm);
		\draw[black] (3,0) circle(1.5cm);
	\end{scope}
	\node[align=center] at (2.3,0) {\textbf{compact case}};
	
	\begin{scope}
		\clip (0,0) ellipse (3cm and 2cm);
		\fill[gray!50!black] (-.3,-2) circle(.5cm);
	\end{scope}
	\node[align=center,text = white] at (-.3,-1.8) {\Large\textbf{A}};
	\begin{scope}
		\clip (0,0) ellipse (3cm and 2cm);
		\fill[gray!50!black] (.7,-2) circle(.4cm);
	\end{scope}
	\node[align=center,text = white] at (.7,-1.8) {\Large\textbf{B}};
		\begin{scope}
		\clip (0,0) ellipse (3cm and 2cm);
		\fill[gray!50!black] (1.7,-1.7) circle(.5cm);
	\end{scope}
	\node[align=center,text = white] at (1.6,-1.5) {\Large\textbf{C}};

\end{tikzpicture}
\caption{Illustrative preview of occurring case distinctions}
\label{fig:operators}
	\end{figure}

An illustration of those operators which are in our focus is given in Figure~\ref{fig:operators} at the end of this introduction. Let us briefly explain here this illustration of operator situations considered in the paper at hand.
If $T : \X \to \Y$ is an injective, bounded linear operator, then $T^*T : \X\to \X$ is positive and self-adjoint and has as a consequence of Halmos' spectral theorem from \cite{h63} a spectral decomposition with a multiplier function $\lambda$ on some measure space $\left(\Omega, \mu\right)$ (see  Section~\ref{sec:Halmos} for details). Under the auspices of the spectral theorem, we will classify the operators $T$ by means of properties of the multiplier function $\lambda$. Our focus is on the classification of {\sl ill-posed} situations that are indicated by the occurrence of {\sl essential zeros} of the function $\lambda$. Such ill-posed problem are characterized in Figure~\ref{fig:operators} by two areas marked hatched in red. The right crescent of this red area is devoted to the subcase of compact operators (type~II ill-posedness), whereas the left part of the red area is reserved for non-compact operators (type~I ill-posedness).   However, we have to exclude the subsituations $(A)$, $(B)$ and $(C)$, marked in black. First the exclave area (A) in black expresses the case of non-compact operators, where the measure $\mu(\Omega)$ is {\sl finite} (Assumption \ref{ass:main1}(a) below is violated).
This case is widely discussed in the literature, in particular for multiplication operators on $\X=L^2(0,1)$ with Lebesgue measure $\mu=\leb$, and we give an overview of literature and specific approaches for that subcase in the appendix. The second black half-circle (B) in Figure~\ref{fig:operators} collects exceptional situations, where the corresponding distribution function $\Phi$ (see formula \eqref{eq:Phi}) is {\sl not informative} and Assumption \ref{ass:main1}(b) below is violated. Such situation is made explicit in Example~\ref{ex:counter2} below. It also occurred in \cite[Example~1]{HofKind10} in form of a diagonal operator in $\ell^2$, where only a subsequence of the diagonal elements tends to zero.  The last small black semi-circle (C) is reserved for possibly occurring situations, where $\mu$ cannot chosen to be the Lebesgue measure $\leb$ on some subset of $\R^d$ (Assumption \ref{ass:main2} is violated). All three cases  marked in black are beyond our focus of this paper.

The paper is organized as follows: In the following Section \ref{sec:compact} we recall the well-known concept of the degree of ill-posedness for  {equations with} compact operators and derive equivalent characterizations of the important \emph{interval of ill-posedness}. {Section~\ref{sec:Halmos} is dedicated to the inclusion of the non-compact case.} We there state the Halmos version of the spectral theorem to be used in the following, define the corresponding distribution function, and illustrate along examples how it can be used to define the degree of ill-posedness. Afterwards, we give our definition and relate it to the compact subcase. Since the focus of our studies is on multiplication operators on \emph{infinite measure spaces}, a glimpse of the finite measure case is given in an appendix. {In Section~\ref{sec:examples} we discuss prominent examples of non-compact operators
with non-closed range,} such as the Hausdorff moment operator and convolutions, and derive their degree of ill-posedness. Section~\ref{sec:unbounded} is briefly devoted to the case of unbounded operators, which in principle could be handled similarly. We end our manuscript by a brief conclusion and outlook in Section~\ref{sec:outlook}.

\section{The degree of ill-posedness for {equations with} compact operators} \label{sec:compact}

To get a blueprint for an application to non-compact operators, we consider in this section \emph{compact} and  injective  linear operators $T:\X \to \Y$ mapping between infinite dimensional Hilbert spaces $\X$ and $\Y$ that lead to ill-posed operator equations \eqref{eq:model} of type~II in the sense of Nashed. They are well-studied in the literature. These studies include reasonable definitions for moderately, severely and mildly ill-posed problems as well as for the degree of ill-posedness, which should be extendable to ill-posendness of type~I.

For injective compact operators $T$ leading to ill-posedness of type~II, the degree of ill-posedness is typically defined by means of the singular value decomposition. Recall, that there exists a uniquely determined \emph{singular system} $\left\{ \left(u_n, v_n, \sigma_n\right)\right\}_{n \in \N}$ of $T$ such that $\left\{u_n\right\}_{n \in \N}$ is a complete orthonormal system in $\X$, $\left\{v_n\right\}_{n \in \N}$ is a complete orthonormal system on the range $\ran(T)$, and $\sigma_1 \geq \sigma_2 \geq ... \geq \sigma_n \geq ...$ are the corresponding singular values obeying the conditions
	\[
	T u_n = \sigma_n v_n, \qquad T^* v_n = \sigma_n u_n \qquad (n \in \N).
	\]
The ill-posedness of {equations involving} $T$ is due to the fact that $\lim_{n \to \infty} \sigma_n = 0$, and the decay rate of $\{\sigma_n\}_{n \in \N}$ indicates the severeness of this ill-posedness.
The first case of a polynomial decay of singular values for compact operators is well-known in the literature as {\sl moderate ill-posedness} (cf., e.g.,~\cite{Hof86})), the second case of an exponential decay as {\sl severe ill-posedness} (cf., e.g.,~\cite{EHN96}). The third case of a logarithmic decay can consequently be characterized as {\sl mild ill-posedness} (cf., e.g,~\cite{AgMa22}).

For a refinement of ill-posedness measurements, in \cite{HofTau97} (see also \cite{HofKind10}) an {\sl interval of ill-posedness} was introduced, which turns out to be  helpful also when generalizing the degree of ill-posedness to non-compact operators. All the above nomenclature can be expressed by means of this interval as collected in the following definition:
\begin{Def}\label{def:degree_compact}
Let $T : \X \to \Y$ be an injective and compact linear operator between infinite-dimensional Hilbert spaces $\X$ and $\Y$, and let $\sigma_n$ be the corresponding singular values. We define the {\bf interval of ill-posedness of {the operator equation \eqref{eq:model}}} as
\begin{equation} \label{eq:interval}
	[A_\sigma,B_\sigma]:=\left[\liminf \limits _{n \to \infty}
	\frac{-\log(\sigma_n)} {\log(n)}\,,\,\limsup \limits _{n \to
		\infty} \frac{-\log(\sigma_n)} {\log(n)}\right]  \subset [0,\infty].
\end{equation}
Furthermore, we call the ill-posed operator equation \eqref{eq:model} \textbf{mildly ill-posed} whenever $A_\sigma = B_\sigma = 0$, \textbf{severely ill-posed} whenever $A_\sigma = B_\sigma = \infty$, \textbf{moderately ill-posed} whenever $0 < A_\sigma \leq B_\sigma < \infty$, and \textbf{ill-posed of degree $s>0$} if $A_\sigma = B_\sigma = s \in (0,\infty)$.
\end{Def}

\begin{ex}[Riemann–Liouville fractional integration of order $\alpha>0$ (cf.~\cite{VuGo94})] \label{ex:frac}
{\rm For the family of the compact operators $T:=J^\alpha: \X=L^2(0,1) \to \Y=L^2(0,1)$ defined as
$$[J^\alpha x](s):= \int_0^{s} \frac{(s-t)^{\alpha-1}}{\Gamma(\alpha)} x(t)\, d t \quad (s \in [0,1]),$$
we have that $\sigma_n \asymp n^{-\alpha}$ as $n \to \infty$, and \eqref{eq:model} is ill-posed of degree $s=\alpha$, since $\lim_{n \to \infty} \frac{-\log(\sigma_n)} {\log(n)}=\alpha$.
}\end{ex}
\begin{ex}[Multivariate integration in the $d$-dimensional case (cf.~\cite{HF23})] \label{ex:multi}
{\rm  For the compact integration operator $T:=J_d: \X=L^2((0,1)^d) \to \Y=L^2((0,1)^d)$ defined as
$$[J_d x](s_1,...,s_d):= \int_0^{s_1} ...\int_0^{s_d} x(t_1,...,t_d)\, d  t_d ... d t_1 \quad ((s_1,...,s_d) \in (0,1)^d),  $$
we derive that $\sigma_n \asymp \frac{[\log(n)]^{d-1}}{n}$ as $n \to \infty$, and \eqref{eq:model} is ill-posed of degree $s=1$
({\sl independent of the dimension} $d \in \N$), because  $\lim _{n \to \infty} \frac{-\log(\frac{[\log(n)]^{d-1}}{n})}{\log(n)}= 1$.
}\end{ex}

The aim of the paper at hand is to extend the above definition to the case of non-compact operators $T$ leading to ill-posed operator equations \eqref{eq:model} of type~I. To do so, several issues have to be considered. First of all, a non-compact operator does not possess a singular value decomposition. This can be handled by considering the self-adjoint operator $T^*T : \X \to \X$ instead, which possesses a spectral decomposition to be discussed below. However, the lack of singular values for $T$ (or of eigenvalues for $T^*T$) still requires further amendments, since no decay rate can be defined directly. The idea we will elaborate therefore can be explained in case of compact linear operators as now described. Following \cite{bhmr07}, we define the \emph{counting function} $\Phi : (0, \infty) \to \R$ as
\begin{equation} \label{eq:Phidisc}
	\Phi(\epsilon)=\#\{n \in \N :\, \sigma_n^2 > \epsilon\} \qquad (\epsilon>0)\,.
\end{equation}
Then the degree of ill-posedness can also be measured by the rate of increase of $\Phi$ as $\epsilon\searrow 0$. Obviously, it holds $\Phi(\epsilon)=0$ for $\epsilon \ge \sigma_1^2=\|T^*T\|$, $0 < \Phi \left(\epsilon\right) < \infty$ for all $0<\epsilon < \sigma_1^2$, and the limit condition $\lim_{n \to \infty} \sigma_n = 0$ indicating ill-posedness transforms into $\lim_{\epsilon \searrow 0} \Phi(\epsilon)=\infty$.

\begin{theorem}
The interval of ill-posedness of {\eqref{eq:model}} is given by
\[
[A_\sigma, B_\sigma] = \left[ \liminf\limits_{\epsilon \searrow 0} \frac{\log(\epsilon)}{-2 \log\left(\Phi(\epsilon)\right)}, \limsup\limits_{\epsilon \searrow 0} \frac{\log(\epsilon)}{-2 \log\left(\Phi(\epsilon)\right)}\right].
\]
\end{theorem}
\begin{proof}
For brevity, we denote
\begin{align*}
A_\Phi & := \liminf_{\epsilon \searrow 0} \frac{\log(\epsilon)}{-2 \log\left(\Phi(\epsilon)\right)},& \qquad B_\Phi & := \limsup_{\epsilon \searrow 0} \frac{\log(\epsilon)}{-2 \log\left(\Phi(\epsilon)\right)}.
\end{align*}
With this notation, we have to show that $[A_\sigma, B_\sigma]  = [A_\Phi,B_\Phi]$, which corresponds to $A_\Phi = A_\sigma$ and $B_\Phi = B_\sigma$.
\begin{description}\item[\fbox{$A_\Phi \geq A_\sigma$}] By definition, for every $\delta > 0$ there exists $N_1 = N_1(\delta)$ such that $ \frac{-\log(\sigma_n)} {\log(n)}  > A_\sigma - \delta$ for all $n \geq N_1$. This implies $\sigma_n < n^{\delta - A_\sigma}$ for all $n \geq N_1$ and hence
	\[
	\left\{n  \in\N \mid n \geq N_1, \sigma_n^2>\epsilon\right\} \subseteq \left\{n  \in\N \mid n \geq N_1, n^{2 \delta-2A_\sigma}>\epsilon\right\}.
	\]
	This implies
	\[
	\Phi \left(\epsilon\right) \lesssim \# \{n \in \N\mid n^{2 \delta-2A_\sigma} > \epsilon\} \asymp \left(\frac1\epsilon\right)^{\frac{1}{2(A_\sigma-\delta)}}
	\]
	as $\epsilon \searrow 0$. Consequently $-2\log(\Phi(\epsilon)) \geq -2\frac{1}{2(A_\sigma-\delta)} \log
	\left(\frac1\epsilon\right) = \frac{1}{A_\sigma-\delta} \log(\epsilon)$ as $\epsilon \searrow 0$, and since $\log(\epsilon) < 0$ in this limiting process, we obtain
	\[
	\liminf\limits_{\epsilon \searrow 0} \frac{\log(\epsilon)}{-2 \log\left(\Phi(\epsilon)\right)} \geq A_\sigma-\delta.
	\]
	Letting $\delta \to 0$ shows $A_\Phi = \liminf_{\epsilon \searrow 0} \frac{\log(\epsilon)}{-2 \log\left(\Phi(\epsilon)\right)} \geq A_\sigma$.
\item[\fbox{$A_\Phi \leq A_\sigma$}] By definition, for every $\delta > 0$ there exists $\epsilon_0 > 0$ such that $\frac{\log(\epsilon)}{-2 \log\left(\Phi(\epsilon)\right)} > A_\Phi-\delta$ for all $0 < \epsilon < \epsilon_0$. This implies (for $\delta < A_\Phi$, which we can assume w.l.o.g.) that $\Phi(\epsilon) < \left(\frac1\epsilon\right)^{\frac{1}{2(A_\Phi-\delta)}}$ for all $0 < \epsilon < \epsilon_0$, and consequently $\sigma_n  < n^{\delta - A_\Phi}$ as $n \to \infty$. This implies $A_\sigma > A_\Phi - \delta$, which by $\delta \searrow 0$ yields $A_\Phi \leq A_\sigma$.
\end{description}
The proof of $B_\Phi = B_\sigma$ follows similarly.  
\end{proof}

\begin{cor}\label{cor:interval}
Let $T : \X \to \Y$ be an injective and compact linear operator between infinite-dimensional Hilbert spaces $\X$ and $\Y$, and let $\Phi$ its counting function as in \eqref{eq:Phidisc}. Then the interval of ill-posedness of {\eqref{eq:model}} is given as
$$
\left[A_\Phi, B_\Phi\right] := \left[ \liminf\limits_{\epsilon \searrow 0} \frac{\log(\epsilon)}{-2 \log\left(\Phi(\epsilon)\right)}, \limsup\limits_{\epsilon \searrow 0} \frac{\log(\epsilon)}{-2 \log\left(\Phi(\epsilon)\right)}\right],
$$
and the operator equation \eqref{eq:model} is mildly ill-posed whenever $A_\Phi = B_\Phi = 0$, severely ill-posed whenever $A_\Phi = B_\Phi = \infty$, moderately ill-posed whenever $0 < A_\Phi \leq B_\Phi < \infty$, and ill-posed of degree $s > 0$ if $A_\Phi = B_\Phi = s$.
 \end{cor}

As mentioned before, the counting function $\Phi$ as introduced by formula \eqref{eq:Phidisc} received some attention when studying statistical inverse problems before, see \cite{bhmr07}. The above corollary is advantageous whenever $\Phi$ can be computed directly without knowing the singular values $\sigma_n$ explicitly. An example for this is the following situation. Consider the negative Laplace operator $-\Delta$ on some compact, smooth Riemannian manifold $S$ of dimension $d$. Then  it follows from Weyl's asymptotic law (see \cite{Taylor}, Chap. 8, Thm. 3.1 and Cor. 3.5) that the counting function
	\[
	N(\lambda) = \# \left\{i \in \N \mid \lambda_i \leq \lambda\right\}
	\] of its eigenvalues $\lambda_i$ satisfies
\[
N(\lambda) \sim \frac{\text{vol}(S)}{\Gamma \left(\frac{d}{2} + 1\right) \left(4 \pi\right)^{\frac{d}{2}}} \lambda^{\frac{d}{2}},
\]
see also Sect. 5 in \cite{bhmr07}. Consequently, if $T^*T$ can be written as $T^*T = \Theta \left(-\Delta\right)$ with a monotonically decreasing function $\Theta$ satisfying $\lim_{t \to \infty} \Theta(t) = 0$, we obtain
\begin{equation} \label{eq:Lapphi}
\Phi\left(\epsilon\right)\sim c_S \left(\Theta^{-1} \left(\epsilon\right)\right)^{\frac{d}{2}}.
\end{equation}

\vspace{-0.4cm}

\begin{ex}[Inverse of negative Laplace operator in the $d$-dimensional case] \label{ex:Laplace}
{\rm  For the special case $\Theta(t)=t^{-2}\;(t>0)$ of formula \eqref{eq:Lapphi}, we find for all $d \in \N$ that the self-adjoint compact operator $T^*T:= \left(-\Delta\right)^{-2}$ on the manifold $S$ corresponds to the counting function $\Phi(\epsilon)\sim c_S\,\epsilon^{-\frac{d}{4}}$.
This yields the limit condition $\lim_{\epsilon \searrow 0} \frac{\log(\epsilon)}{-2\log(\Phi(\epsilon))}=\frac{2}{d}$ and indicates that the problem is ill-posed of degree $s=2/d$, which (in contrast to the multidimensonal Example~\ref{ex:multi}) {\sl strongly depends on the space dimension} $d$ of the manifold $S$ under consideration.
}\end{ex}

The above theorem and the derived corollary show that we can express the degree of ill-posedness in the compact case equivalently by means of the counting function $\Phi$ from \eqref{eq:Phidisc}. In the following, we will use this as a blueprint to consistently extend the corresponding definition to the non-compact case under certain restrictions.

\begin{rem} {\rm For Example~\ref{ex:Laplace}, solving the operator equation \eqref{ass:main1} requires to find $u=(-\Delta)g$, which means that {\sl second} partial derivatives of $g$ have to be calculated. The factor $2$ in the degree of ill-posedness $s=2/d$ seems to express this. Namely, for
the compact embedding operator $T:H^p((0,1)^d) \to L^2((0,1)^d)$  from the Sobolev space of degree $p>0$ to $L^2$ in $d$ dimensions we have singular values as $\sigma_n \asymp n^{-\frac{p}{d}}$ (cf.,e.g.,\cite[\S 3c]{Koenig86}) implying a degree of ill-posedness as $s=p/d$, and just for $p=2$ calculations up to second partial derivatives occur.}
\end{rem}

\section{The spectral theorem as a tool for defining the degree of ill-posedness also in the non-compact case}\label{sec:Halmos}

At this point applying the adjoint operator $T^*: \Y \to \X$ for $T$ to the operator equation \eqref{eq:model}, we restrict our considerations to the occurring normal equation  $T^*Tu=h$ for $h=T^*g$ to \eqref{eq:model} with the bounded, self-adjoint and positive semi-definite operator $T^*T: \X \to \X$.
Let us mention that the ill-posedness and the type of ill-posedness of equation \eqref{eq:model} carry over to the normal equation.

\subsection{The spectral theorem according to Halmos}

Moreover, we recall the spectral theorem in multiplication operator form  by Halmos \cite{h63} and take into account the additional notes in \cite[p.~47]{Haase18} and \cite[Chap.~VII.1]{w05}. Then the following holds true for the normal operator $T^*T: \X \to \X$: There exists a locally compact topological space $\Omega$, a semi-finite (Radon) measure $\mu$ on $\Omega$ with finite or infinite measure value $\mu(\Omega)$,
a real-valued function $\lambda \in L^\infty\left(\Omega, \mu\right)$,
and a unitary mapping $W : L^2 \left(\Omega, \mu\right) \to \X$ such that
\begin{equation} \label{eq:Multop}
	W^*T^*TW = M_\lambda\,,
\end{equation}
with the multiplication operator $M_\lambda :L^2 \left(\Omega, \mu\right) \to L^2 \left(\Omega, \mu\right)$ defined for all $\xi \in L^2(\Omega,\mu)$ as
\begin{equation} \label{eq:lambda}
	[M_\lambda \xi](\omega):=\lambda(\omega)\cdot \xi(\omega) \qquad (\mu-\mbox{a.e.~on}\;\Omega).
\end{equation}
We note that this existence assertion applies to \emph{compact and non-compact operators} $T:\X \to \Y$ in the same way.
Due to the positive semi-definiteness of the operator $T^*T: \X \to \X$, the multiplier function $\lambda$ is non-negative almost everywhere. Furthermore, its \emph{essential range}
$${\rm essran}(\lambda)=\{z \in [0,\infty) :\, \mu(\{\omega \in \Omega:\, |z-\lambda(\omega)| \le \eta\})>0 \;\mbox{for all}\; \eta>0\}$$
coincides with the spectrum ${\rm spec}(T^*T)$ of the operator $T^*T$ (cf., e.g., \cite[Theorem~2.1(g)]{Haase18}) and is hence a subset of $[0,\|T^*T\|]=[0,\|T\|^2]$.
Let us mention that the injectivity requirement for $T$ made throughout this paper implies that $\mu(\{\omega \in \Omega:\,\lambda(\omega)=0\})=0$.

{The operator equation \eqref{eq:model} is ill-posed if and only if $0 \in {\rm spec}(T^*T)$. Then as a consequence of the above mentioned coincidence between spectrum of $T^*T$ and the essential range of the multiplier function $\lambda$, the criterion of ill-posedness for \eqref{eq:model}
is the occurrence of an \emph{essential zero} of the function $\lambda$, which means that zero belongs to the essential range of $\lambda$.} Our goal is to characterize the ill-posedness of the original operator equation \eqref{eq:model} and its associated normal equation by properties of this multiplier function $\lambda$ in \eqref{eq:lambda} and their implications. However, this is not immediately possible without further restrictions, since the spectral decomposition \eqref{eq:Multop} is not uniquely determined for the operator $T$ as the following lemma indicates. Note that this observation was already employed in \cite{bhmr07}:
\begin{lem}\label{lem:non_uniqueness}
	If $\kappa : \Omega \to (0,\infty)$ is a $\mu$-measurable function, then also
	\[
	\tilde W^*T^*T\tilde W = M_\lambda\,,
	\]
	with
	\[
	\tilde \mu := \kappa \mu, \qquad \tilde W :  L^2 \left(\Omega, \tilde \mu\right) \to \X, \quad \tilde W = W \circ M_{\kappa^{\frac12}}
	\]
	is a spectral decomposition of $T^*T$.
\end{lem}
\begin{proof}
	First we show that $\tilde W$ is unitary. For $f \in L^2 \left(\Omega, \tilde \mu\right)$ we have
\begin{multline*}
	\left\langle \tilde W f, \tilde W g\right\rangle_{\X} = \left\langle W (\kappa^{\frac12} f), W (\kappa^{\frac12} g)\right\rangle_{\X} = \left\langle \kappa^{\frac12} f,\kappa^{\frac12} g \right\rangle_{L^2(\Omega,\mu)} \\= \int_{\Omega} \kappa f \bar g \,\mathrm d \mu = \int_{\Omega} f \bar g \,\mathrm d \tilde \mu = \left\langle f,g\right\rangle_{L^2 \left(\Omega, \tilde \mu\right)}
\end{multline*}
	since $\kappa$ is real-valued and $W$ unitary. Furthermore, the computation
	\begin{multline*}
		\left\langle \tilde W f,x\right\rangle_{\X} = \left\langle W \left(\kappa^{\frac12} f\right),x\right\rangle_{\X} = \left\langle \kappa^{\frac12} f, W^* x \right\rangle_{L^2 \left(\Omega, \mu\right)}\\
		= \int_{\Omega}  \kappa^{\frac12} f (W^* x)\,\mathrm d \mu = \int_{\Omega}  \kappa^{-\frac12} f (W^* x)\,\mathrm d \tilde \mu = \left\langle  f, \kappa^{-\frac12}W^* x \right\rangle_{L^2 \left(\Omega, \tilde \mu\right)}
	\end{multline*}
	shows that the adjoint of $\tilde W$ is given by $\tilde W^* = M_{\kappa^{-\frac12} }\circ W^* $. This yields finally
	\[
	\tilde W^*T^*T\tilde W = M_{\kappa^{-\frac12}} W^*   T^*T W M_{\kappa^{\frac12}} = M_{\kappa^{-\frac12}} M_\lambda M_{\kappa^{\frac12}} = M_\lambda.
	\]
\end{proof}

\begin{rem}{\rm
The above non-uniqueness will lead to severe problems when defining the degree of ill-posedness for {equations with} non-compact $T$, {but also in the compact case it brings in a subjective factor}. To illustrate this, recall that our definition of the degree of ill-posedness for compact operators via the counting function $\Phi$ from \eqref{eq:Phidisc} explicitly relies on the {\sl counting measure} $\mu_{disc}=\#$ in \eqref{eq:Phidisc}. If we allow for other measures, we will obtain completely different asymptotical behaviors of $\Phi$. As an example we next consider a history match problem.}
\end{rem}

\begin{ex}[Periodic backwards heat equation]\label{ex:bhe}
{\rm
\begin{align*}
	\frac{\partial u}{\partial t} u(x,t) & = \frac{\partial^2}{\partial t^2} \left(x,t\right) \qquad &\text{in } &(-\pi,\pi] \times (0,\bar t),\\
	u(x,0) & = f(x) \qquad & \text{on } &[-\pi,\pi],\\
	u(-\pi,t) & = u(\pi,t)\qquad &\text{on }& (0,\bar t],
\end{align*}
with a given terminal time $\bar t > 0$, and denote by $T_{\mathrm{BHE}} : L^2 \left({-\pi},\pi\right) \to L^2 \left(-\pi,\pi\right)$. The corresponding operator equation \eqref{eq:model} is widely understood to be exponentially ill-posed, and this is reflected by the series representation
\[
(T_{\mathrm{BHE}}f)(x) = \sum_{k \in \Z} \exp\left(-k^2 \bar t\right) \exp\left(\textup{i} kx\right) \hat f(k)\qquad(x \in [-\pi,\pi]),
\]
with the Fourier coefficients $\hat f(k)$ of $f$. Note, that this representation corresponds to the spectral decomposition
\[
\mathcal F^* T_{\mathrm{BHE}} \mathcal F = M_\lambda
\]
with $\lambda \in \ell^2(\mathbb Z)$ given by $\lambda(k) = \exp\left(-k^2 \bar t\right)$. The corresponding distribution function $\Phi$ from \eqref{eq:Phidisc} is given by
\[
\Phi(\epsilon) = \# \{k \in \Z \mid \lambda(k) > \epsilon\} \asymp 2 \sqrt{\frac{1}{\bar t} \log \left(\frac{1}{\epsilon}\right)} \qquad\text{as}\qquad \epsilon \searrow 0.
\]
However, if we choose $\kappa(k) := \exp\left(\bar t k^2\right)$ in Lemma \ref{lem:non_uniqueness}, then we obtain a different measure $\tilde \mu = \kappa \mu$ (actually a weighted counting measure) such that
\[
\tilde \Phi(\epsilon) = \tilde \mu \left(\left\{k \in \Z \mid \lambda(k) > \epsilon \right\}\right) = \sum_{k \in\Z, k^2 \leq \frac1{\bar t} \log\left(\frac1\epsilon\right)} \kappa(k) \sim \frac{1}{\epsilon}.
\]
This indicates a completely different degree of ill-posedness, and hence the counting measure has to be fixed {as benchmark measure for compact operators} to obtain a meaningful concept.
}\end{ex}

\begin{rem}
{\rm Note, that no further assumptions on the function $\kappa$ in Lemma \ref{lem:non_uniqueness} are required, i.e. for any measurable $\kappa$, the corresponding $\tilde \mu = \kappa \mu$ is again a measure. Furthermore, even the additional properties of the measure space $\left(\Omega,\mu\right)$ gained from the spectral theorem, i.e. semi-finiteness, can be transferred to $\tilde \mu$ by requiring that $\kappa$ is continuous. In contrast, by specific choices of $\kappa$ we can obtain $\tilde \mu(\Omega) < \infty$ even for $\mu(\Omega) = \infty$ and vice versa. In Figure \ref{fig:operators} we have indicated that we will exclude the case $\mu(\Omega) < \infty$, and hence the above consideration again indicates that we have to fix the measure $\mu$ later on in order to derive a meaningful concept of the degree of ill-posedness.
}\end{rem}

\subsection{The distribution function measuring the strength of ill-posedness}

Inspired by the counting function \eqref{eq:Phidisc} we will consider, as an analog under the auspices of the spectral theorem and now applicable to both compact and non-compact operators $T$,  the non-negative, non-increasing and right-continuous \emph{distribution function}
\begin{equation}\label{eq:Phi}
	\Phi_{\lambda,\mu}(\epsilon) := \mu \left( \{\omega \in \Omega :\, \lambda(\omega) > \epsilon\}\right) \quad (\epsilon >0).
\end{equation}
This function $\Phi_{\lambda,\mu}$ is not in general informative, since it cannot be excluded that it might attain the value $+\infty$ on $(0,\infty)$:
\begin{ex}[Non-informative distribution function] \label{ex:counter2}
	{\rm Consider $\Omega=[0,\infty)$ and the associated  Lebesgue measure $\mu = \leb$ on $\R$ with $\leb(\Omega)=\infty$. Then the function $\lambda(\omega)=\sin^2(\omega) \;(\omega \in [0,\infty))$ leads to
$$\Phi_{\lambda,\leb}(\epsilon)=\left\{\begin{array}{ccc} \infty & \mbox{for} & 0<\epsilon <  1,\\ 0 & \mbox{for} & \epsilon \ge 1 \end{array}\right.$$
as distribution function. Hence $\Phi_{\lambda,\leb}$ is not informative at all. This is, however, an interesting example, because $0 \in {\rm essran}(\lambda) = {\rm spec}(T^*T)$ indicates ill-posedness.}
\end{ex}

As a consequence, we have to pose assumptions on the measure space $(\Omega,\mu)$ and the multiplier function  $\lambda$ such that the function $\Phi_{\lambda,\mu}$ is finite and hence informative. Furthermore, note that the function $\Phi_{\lambda,\mu}$ is determined by both $\lambda$ and $\mu$, and not solely by the operator $T$. The induced spectral decompositions in Lemma~\ref{lem:non_uniqueness} might yield different functions $\Phi_{\lambda,\mu}$ for the same operator $T$. This corresponds to the problem outlined in Example~\ref{ex:bhe} and indicates that we need to restrict to a certain \textit{benchmark measure space} $\left(\Omega, \mu\right)$.

\begin{ass} \label{ass:main1}
	Let, for the operator $T^*T$, for the measure space $(\Omega,\mu)$  and for the multiplier function $\lambda$ in \eqref{eq:lambda} defining the multiplication operator $M_\lambda$ from \eqref{eq:Multop} obey the following two assumptions:
	\begin{itemize} 
		\item[(a)] $\mu(\Omega)=\infty$.
		\item[(b)] For the distribution function $\Phi_{\lambda,\mu}$ it holds $0 \le \Phi_{\lambda,\mu} \left(\epsilon\right) < \infty$ for all $\epsilon > 0$ with the limit condition  $\lim \limits_{\epsilon \searrow 0} \Phi_{\lambda,\mu}(\epsilon)=\infty$.
	\end{itemize}
\end{ass}

\begin{pro}[Ill-posedness] \label{pro:illposed}
	If Assumption~\ref{ass:main1} is satisfied for the operator $T$ and its corresponding measure space $(\Omega,\mu)$ and multiplier function $\lambda$, then the operator equation \eqref{eq:model} is ill-posed.
\end{pro}
\begin{proof}
	We show that $0  \in {\rm essran}(\lambda)$, which is equivalent to  $0\in  {\rm spec}(T^*T)$ and indicates the ill-posedness of \eqref{eq:model}. If $\mu(\{\omega \in \Omega:\,\lambda(\omega) \le \eta\})=0$ were valid for any $\eta>0$ then, in view of $\mu(\Omega)=\infty$ and in contradiction to item (b) of  Assumption~\ref{ass:main1},
	$\Phi_{\lambda,\mu}(\eta)=\mu(\{\omega \in \Omega:\,\lambda(\omega) > \eta\})=\infty$ would hold. Thus we have $\mu(\{\omega \in \Omega:\,\lambda(\omega) \le \eta\})>0$ for all $\eta>0$, which means that $0  \in {\rm essran}(\lambda)$.
\end{proof}

Let us briefly discuss sufficient conditions and counterexamples for item (b) in Assumption~\ref{ass:main1}.

\begin{lem} \label{lem:ffinite}
	Let $\mu(\Omega)=\infty$. Then the following condition is sufficient for obtaining the assertion of item (b) in Assumption~\ref{ass:main1}:
	\begin{itemize} 
		\item[(c)] There exists a $\mu$-measurable, non-negative and non-decreasing function $f(\zeta)\;(0 \le \zeta < \infty)$ with $f(\zeta)>0 $ for $\zeta>0$ such that it holds $\int \limits_\Omega f(\lambda(\omega)) d\mu <\infty$.
	\end{itemize}
\end{lem}
\begin{proof}
	For $\epsilon>0$ we can estimate due to the monotonicity of $f$ and $f(\epsilon)>0$  as follows:
	\begin{multline*}0 \le \Phi_{\lambda,\mu}\left(\epsilon\right) =\int \limits_{\{\omega \in \Omega:\,\epsilon < \lambda(\omega)\}} 1 \d \mu= \frac{1}{f(\epsilon)} \int\limits_{\{\omega \in \Omega:\,\epsilon < \lambda(\omega)\}} f(\epsilon) \d \mu \\\leq  \frac{1}{f(\epsilon)} \int\limits_{\{\omega \in \Omega:\,\epsilon < \lambda(\omega)\}} f(\lambda(\omega)) \d \mu\leq\quad  \frac{1}{f(\epsilon)} \int \limits_\Omega f(\lambda(\omega)) \d \mu  < \infty.
	\end{multline*}
	Note that the last inequality is a consequence of the nonnegativity of the function $\lambda$ and yields the finiteness of $\Phi_{\lambda,\mu}(\epsilon)$ for $\epsilon>0$. The limit condition arises from
$\lim _{\epsilon \searrow 0}\Phi_{\lambda,\mu}(\epsilon)=\mu(\Omega)=\infty$ if we let $\epsilon$ tend to zero
	in \eqref{eq:Phi}.
\end{proof}

An important example for condition (c) in Lemma~\ref{lem:ffinite} is the setting $f(\zeta):=\zeta^p$ for exponents \linebreak $0 < p < \infty$, for which Lemma~\ref{lem:ffinite} under $\mu(\Omega)=\infty$ ensures the validity of Assumption~\ref{ass:main1}. In particular, this is the case
for $\lambda \in L^p(\Omega,\mu)\;(1 \le p<\infty)$, and as a corollary to Lemma~\ref{lem:ffinite} we immediately find from this lemma:

\begin{cor}[Multiplier functions of $L^p$-type] \label{cor:Lp}
	If we have, under the condition $\mu(\Omega)=\infty$, that $\lambda \in L^p(\Omega,\mu)$ is true for some $1 \le p < \infty$, then item (b) of Assumption~\ref{ass:main1} is satisfied.
\end{cor}

Note that for general operators $T$ the assumption $\lambda \in L^p \left(\Omega, \mu\right)$ for some $p\in [1,\infty)$ seems reasonable. In Lemma 2 in \cite{bhmr07} it is shown that one can choose $\lambda \in L^1 \left(\Omega,\mu\right)$ as soon as $T$ is injective and $T^*$ transforms white noise into a Hilbert space valued random variable. Finally, allowing for $p > 1$, the assumption $\lambda \in L^p(\Omega,\mu)$ generalizes this to more regular, colored noise.

\medskip

Example \ref{ex:counter2} above and the following additional counterexample show that the assertion of item~(b) in Assumption~\ref{ass:main1} may be violated for multiplier functions $\lambda$.

\begin{ex} \label{ex:counter1}
	{\rm Let $\mu(\Omega)=\infty$ and $\lambda(\omega)=c>0$ for all $\omega \in \Omega$. Then we have for the distribution function $\Phi_{\lambda,\mu}(\epsilon)=\left\{\begin{array}{ccc} \infty & \mbox{for} & 0<\epsilon < c\\ 0 & \mbox{for} & \epsilon \ge c \end{array}\right.,$ and hence item~(b) in Assumption~\ref{ass:main1} is violated, but we have no ill-posedness, because $0 \notin {\rm essran}(\lambda)$.}
\end{ex}

In case that Assumption~\ref{ass:main1} holds true, it makes sense to formulate an analog to Definition~\ref{def:degree_compact} under the auspices of the spectral theorem as follows:
\begin{Def}\label{def:degree_noncompact}
Let $\X$ and $\Y$ denote infinite-dimensional Hilbert spaces. Moreover, let $T : \X \to \Y$ be an injective and bounded linear operator with non-closed range, each mapping between $\X$ and $\Y$ such that the Assumption \ref{ass:main1} is satisfied w.r.t.~the distribution function $\Phi_{\lambda,\mu}$ from  \eqref{eq:Phi}. Then the \textbf{interval of ill-posedness of {the operator equation \eqref{eq:model}}} is defined as
\begin{equation} \label{eq:interval2}
\left[A_{\Phi_{\lambda,\mu}}, B_{\Phi_{\lambda,\mu}}\right] := \left[ \liminf\limits_{\epsilon \searrow 0} \frac{\log(\epsilon)}{-2 \log\left(\Phi_{\lambda,\mu}(\epsilon)\right)}, \limsup\limits_{\epsilon \searrow 0} \frac{\log(\epsilon)}{-2 \log\left(\Phi_{\lambda,\mu}(\epsilon)\right)}\right],
\end{equation}
and the operator equation \eqref{eq:model} is called \textbf{mildly ill-posed} whenever $A_{\Phi_{\lambda,\mu}} = B_{\Phi_{\lambda,\mu}} = 0$, \textbf{severely ill-posed} whenever $A_{\Phi_{\lambda,\mu}} = B_{\Phi_{\lambda,\mu}} = \infty$, \textbf{moderately ill-posed} whenever $0 < A_{\Phi_{\lambda,\mu}} \leq B_{\Phi_{\lambda,\mu}} < \infty$, and \textbf{ill-posed of degree $s > 0$} if $A_{\Phi_{\lambda,\mu}} = B_{\Phi_{\lambda,\mu}} = s$. Then the distribution functions
\end{Def}
\begin{theorem}
For compact operators $T$, the intervals and degrees of ill-posedness from Definitions~1 and 2 agree.
\end{theorem}
\begin{proof}
For compact operators $T$ with singular values $\{\sigma_n\}_{n=1}^\infty$ there is a unitary operator $W: L^2(\N,\mu_{disc}) \to \X$ for $\Omega=\N$ and the counting measure $\mu=\mu_{disc}$ such that the formulas \eqref{eq:Multop} and \eqref{eq:lambda} hold with the multiplier function $\lambda(n)=\sigma_n^2\;(n \in \N)$.
{Then the counting function $\Phi$ from \eqref{eq:Phidisc} and the distribution function $\Phi_{\lambda,\mu}$ from  \eqref{eq:Phi} match. Taking into account the result of Corollary~\ref{cor:interval}, this proves the theorem.}
\end{proof}

Indeed, by means of the spectral decomposition \eqref{eq:Multop}, all information of the equation $T^*Tu=h$ with self-adjoint operator $T^*T$ mapping in the Hilbert space $\X$, which appears as normal equation to \eqref{eq:model}, is contained in the transformed multiplication operator equation
\begin{equation}\label{eq:trans}
	\lambda(\omega)\cdot \xi(\omega)= \zeta(\omega) \qquad (\mu-\mbox{a.e.~on}\;\Omega)
\end{equation}
in $\Omega$, because the functions $\xi,\zeta \in L^2(\Omega,\mu)$ appear by unitary transformation from $u$ and $h$. Formally, the inversion of equation \eqref{eq:trans} requires the division $\xi(\omega)=\frac{\zeta(\omega)}{\lambda(\omega)}\;(\omega \in \Omega, \;\mu-\mbox{a.e.})$. Hence,
instability grows if the decay rates to zero of the multiplier function $\lambda$ become more and more pronounced. This motivates Definition~\ref{def:degree_noncompact} {for a classification of ill-posedness strength within the class of compact operators under the measure space $(\N,\mu_{disc})$ as well as for a classification within families of non-compact operators under the measure space $(\R,\mu_{Leb})$, which will be discussed for example in the following subsection.}

\subsection{Ill-posedness of {equations with} multiplication operators on the real axis} \label{sec:realaxis}

{To enter the world of ill-posed linear operator equations \eqref{eq:model} for injective, bounded and \emph{non-compact} operators $T:\X \to \Y$ with non-closed range,} we first consider in this subsection {\sl pure multiplication operators} $T:=M_\lambda$ defined as
 $[M_\lambda\, \xi](\omega):=\lambda(\omega)\cdot \xi(\omega)\;(\mu-\mbox{a.e.~on}\;\Omega)$ for $\X=\Y=L^2(\Omega)$ on the real axis with a measure space $(\Omega,\leb)$, where $\Omega = \R=(-\infty,+\infty)$ is combined with the Lebesgue measure $\leb$ over $\R$.
We assume that the multiplier functions $\lambda$ is continuous everywhere on $\R$ and obeys the limit condition $\lim_{\omega \to \pm\infty} \lambda(\omega)=0$. Then we always have a non-compact operator $T^*T$ with purely continuous spectrum
${\rm spec}(T^*T)=[0,\|T^*T\|]={\rm essran}(\lambda)=[0,\max \limits_{\omega \in \R}\lambda(\omega)]$ that implies ill-posedness of type~I in the sense of Nashed of the associated operator equation \eqref{eq:model}.
Definition~\ref{def:degree_noncompact} now applies here for this class of non-compact operators with non-closed range by considering the interval of ill-posedness \eqref{eq:interval2}
derived from the distribution function \eqref{eq:Phi} with $\mu=\mu_{Leb}$. As the variants in Example~\ref{ex:Frank} below will show, also here moderately, severely and mildly ill-posed problems can be distinguished intuitively.

\bigskip

\begin{ex} \label{ex:Frank}
	{\rm \begin{itemize} \item[]
			\item[(a)] {\sl Moderate ill-posedness:}
			
			\smallskip
			
			(a1) $\; \lambda (\omega) = \frac{1}{(1+\omega^2)^s} \quad(\omega \in \R,\; s > 0)\quad \mbox{implies}\quad \Phi_{\lambda,\leb}(\epsilon) = 2\mu(\{\omega \ge 0:\,1+\omega^2 \le \epsilon^{-1/s}\})\,,$ hence
			$\; \Phi_{\lambda,\mu_{Leb}}(\epsilon) = 2 \sqrt{\epsilon^{-\frac1s} -1} \asymp \epsilon^{-\frac{1}{2s}}$ as $\epsilon \to 0$. The corresponding interval \eqref{eq:interval2} hence reduces to a singleton  $\left[A_{\lambda,\leb},B_{\lambda,\leb}\right] = \{s\}$, indicating a moderate degree $s>0$ of ill-posedness for the underlying operator equation here.
			
			\smallskip
			
			(a2) $\; \lambda (\omega) = \frac{\omega^2}{1+\omega^4} \quad(\omega \in \R,\;s=1)\quad \mbox{implies} \quad \Phi_{\lambda,\leb}(\epsilon) = 2\mu(\{\omega > 0:\,\omega^2+\omega^{-2} \le \epsilon^{-1}\})\,$, but for $\epsilon\ll1$ we have
			$$\Phi_{\lambda,\leb}(\epsilon) \approx  2\left(\epsilon^{-\frac{1}{2}}-\epsilon^{\frac{1}{2}}\right),\mbox{where the zero of}\; \lambda(\omega)\,\mbox{at}\;\omega=0\;\, \mbox{causes the negative term}. $$
			The impact of this inner zero on the asymptotics of $\Phi_{\lambda,\leb}(\epsilon)$ as $\epsilon \to 0$ is negligible (cf. also \cite[Prop.~1]{mnh22}), and due to the zeros of $\lambda(\omega)$ at infinity $\omega \to \pm \infty$
			we have $\Phi_{\lambda,\leb}(\epsilon) \asymp \epsilon^{-\frac{1}{2}}$ as well as
			$\left[A_{\lambda,\leb},B_{\lambda,\leb}\right] = \{1\}$, expressing a degree one of ill-posedness.
			
			\item[(b)] {\sl Severe ill-posedness:}
			
			\smallskip
			
			$\lambda (\omega) =\exp\left(-|\omega|^s\right) \quad(\omega \in \R,\; s > 0)\quad \mbox{implies}\quad \Phi_{\lambda,\leb}(\epsilon) = 2\mu(\{\omega \ge 0:\,\omega^{s} \le -\log(\epsilon)\})\;$ and
			$ \Phi_{\lambda,\leb}(\epsilon) \asymp [\log(1/\epsilon)]^\frac{1}{s}$ as $\epsilon \to 0\;$
            The zeros of $\lambda(\omega)$ at infinity $\omega \to \pm\infty$ imply $\left[A_{\lambda,\leb},B_{\lambda,\leb}\right] = \{\infty\}$, which is in agreement with the intuitive severe ill-posedness of the underlying operator equation.
			
			\bigskip
			
			\item[(c)] {\sl Mild ill-posedness:}
			
			\smallskip

			$\lambda (\omega) =\left\{\begin{array}{ccl} [\log(|\omega|)]^{-2s} & \mbox{for} & |\omega| \ge e\\1 & \mbox{for} &  |\omega| < e \end{array}\right.\quad (s > 0)\qquad \mbox{implies for} \;\;  0<\epsilon\ll 1$
			$$ \Phi_{\lambda,\leb}(\epsilon) = 2\mu(\{\omega \ge 0:\,\omega^{2s} \le -\log(\epsilon)\})\;(\epsilon\ll 1) \qquad \mbox{and hence}$$
			$\Phi_{\lambda,\leb}(\epsilon) \asymp \exp\left[\left(\frac{1}{\epsilon}\right)^{\frac{1}{2s}}\right]$ as $\epsilon \to 0$, i.e. the interval in \eqref{eq:interval2} is given by $\left[A_{\lambda,\leb},B_{\lambda,\leb}\right] = \{0\}$. Again the zeros of $\lambda(\omega)$ at infinity $\omega \to \pm\infty$, which are of logarithmic type, determine the intuitive mild degree of ill-posedness.

		\end{itemize}
		
		\bigskip
		
		{\parindent0em For} all $s>0$, the variants (a1), (b) and (c) of Example~\ref{ex:Frank}  satisfy Assumption~\ref{ass:main1}, and so does variant (a2). Furthermore, the multiplier function $\lambda$ belongs to $L^p(\Omega,\mu)$ in the sense of Corollary~\ref{cor:Lp}, for the moderately ill-posed case (a)
		whenever $s>\frac{1}{2p}$. It belongs for all $s>0$ even to $L^1(\Omega,\mu)$ in the severely ill-posed case (b). For the mildly ill-posed case (c), it makes sense to the recall Lemma~\ref{lem:ffinite} and in particular the condition (c) ibid. Then the increasing function $f(\zeta)=\exp\left( -2/\zeta^{\frac{1}{2s}}\right)$ is an appropriate choice for $\lambda(\omega)=[\log(\omega)]^{-2s}\;(\omega \ge e)$ from (c) such that Lemma~\ref{lem:ffinite} applies, because $\int \limits_{\R} f(\lambda(\omega)) d\omega < \infty$. In general, such function $f$ can be found if
		a setting $f(\zeta):=1/\lambda^{-1}(\zeta)$ for sufficiently large $\zeta>0$ is useful, which requires a strict decay of $\lambda(\omega)$ as $\omega$ tends to infinity.
	} \hfill\fbox{}\end{ex}

\subsection{{Chances and limitations of Lebesgue measure spaces}} \label{sec:Leb}

{Lemma~\ref{lem:non_uniqueness} has shown that the strength of ill-posedness expressed by interval and degree of ill-posedness in the sense of Definitions~\ref{def:degree_compact} and \ref{def:degree_noncompact} strongly depends on the choice of the measure space $(\Omega,\mu)$.
In order to overcome the difficulties arising from the ambiguity of possible measure spaces, we have great expectations of \emph{benchmark measure spaces}. By using such benchmarks, the outcome of Definition~\ref{def:degree_noncompact} is uniquely determined for the operator $T$.  The use of $(\N,\mu_{disc})$ with the counting measure $\mu_{disc}$ for the classification of compact operators is fairly uncontroversial.
Likewise, the examples in Subsection~\ref{sec:realaxis} establish the Lebesgue measure $\mu_{Leb}$ on $\R$ as a benchmark measure for the classification of families of non-compact operators, which motivates the following Assumption~\ref{ass:main2} as a supplement to  Assumption~\ref{ass:main1}.
\begin{ass}\label{ass:main2}
	Let $T : \X \to Y$ be a injective bounded non-compact linear operator with non-closed range mapping between infinite-dimensional Hilbert spaces $\X$ and $\Y$ such that there exists a spectral decomposition \eqref{eq:Multop} with a convex  set $\Omega \subseteq \R^d$, $\mu(\Omega)=\infty$ and $\mu = \leb$.
\end{ass}
This assumption might be restrictive at the moment, but we will discuss below in Section~\ref{sec:examples} several examples with non-compact operators satisfying this assumption, e.g. convolution operators and the Hausdorff moment problem. Let us emphasize that the Lebesgue measure employed in Assumption~\ref{ass:main2} is natural. Taking note of the transformations as in Lemma \ref{lem:non_uniqueness}, the restriction to a situation with Lebesgue measure has the rationale background that its translation invariance ensures all frequencies to have the same influence.}

{It would be desirable to find a common measure space $(\Omega,\mu)$ which allows for a comparison of the ill-posedness strength \emph{between} compact and non-compact operator.
For example, an open question seems to be whether a severely ill-posed non-compact operator exists which in any sense proves to be more ill-posed than another moderately ill-posed compact operator.
Trying to extend Assumption~\ref{ass:main2} to compact cases, it had been an idea to use $\Omega=[0,\infty)$ with the Lebesgue measure $\mu_{Leb}$ over $\R$ as possible common measure space by exploiting the (generalized) inverse
\begin{equation} \label{eq:invphi}
			\lambda^*(\omega):=\inf\left\{\tau>0:\,\Phi_{\lambda,\mu}(\tau) \le \omega\right\}\quad (0 \le \omega <\infty)
\end{equation}
of the distribution function $\Phi_{\lambda,\mu}$. This \emph{decreasing rearrangement} of the multiplier function $\lambda$ is defined on this subset $[0,\infty)$ of the real axis, even if $\lambda$ is originally defined on a subset of $\R^d$ for $d>1$.
Moreover, $\lambda$ and $\lambda^*$ are closely related, because they have the same distribution function $\Phi_{\lambda,\mu}(\epsilon)$, and the growth rate of this distribution functions as $\epsilon \searrow 0$, which is relevant for Definition~\ref{def:degree_noncompact},
can be rewritten under Assumption~\ref{ass:main1} as decay rate of $\lambda^*(t)$ as $t \to \infty$, see in this context also \cite{mnh22}.}

{Following this idea, compact operators $T=T_{comp}$ with singular values $\{\sigma_n\}_{n=1}^\infty$ and distribution (counting) function $\Phi$ from \eqref{eq:Phidisc} possess the non-negative, non-increasing and  piece-wise constant multiplier function $\lambda_{comp}$ defined by formula
\begin{equation} \label{eq:step}
\lambda_{comp}(\omega)=\{\sigma_n^2 \;\mbox{if}\; n-1 \le \omega < n \;(n=1,2,...)\} \quad (0 \le \omega <\infty)
\end{equation}
as inverse of $\Phi(\epsilon)$ tending to zero as $\omega \to \infty$. This function from \eqref{eq:step} is the decreasing rearrangement of the original multiplier function $\lambda_{comp}(n)=\sigma_n^2\;(n \in \N)$ assigned to $T_{comp}^*T_{comp}$. By bringing the Lebesgue measure over $\R$ into play, we find that
$$\Phi_{\lambda_{comp},\mu_{Leb}}(\epsilon)= \leb(\{\omega \in [0,\infty):\,\lambda_{comp}(\omega)>\epsilon\}) \quad (\epsilon>0)$$
coincides with the corresponding counting function~\eqref{eq:Phidisc}.
 If we consider the corresponding multiplication operator $M_\lambda$ mapping in $L^2([0,\infty),\mu_{Leb})$ defined by formula \eqref{eq:lambda} for the multiplier function $\lambda=\lambda_{comp}$ from \eqref{eq:step}, then this non-negative self-adjoint multiplication operator $M_{\lambda_{comp}}$ has a pure
point spectrum, and the eigenvalues of $M_{\lambda_{comp}}$ are the same values as those of the compact operator $T_{comp}^*T_{comp}$. Namely, the essential range of the multiplier function
and the spectrum values for a multiplication operator are identical. However, the following general assertion of Proposition~\ref{pro:nocompact} shows the limitation of the Lebesgue measure in this setting. Here, $M_{\lambda_{comp}}$ is a non-compact operator with infinite dimensional
eigenspaces for all eigenvalues, and the spectrum of $M_{\lambda_{comp}}$ is not a discrete spectrum with only finite dimensional eigenspaces, which would be necessary for being a compact self-adjoint operator like $T_{comp}^*T_{comp}$.  As a consequence we have that  $M_{\lambda_{comp}}$
as a non-compact operator  cannot be derived from the compact operator  $T_{comp}^*T_{comp}$ by a unitary transformation, This implies that a desired comparison of compact and non-compact operators fails on this way.}

{
\begin{pro} \label{pro:nocompact}
There is \emph{no compact} injective non-negative multiplication operator
$$M_\lambda: L^2([0,\infty),\mu_{Leb}) \to L^2([0,\infty),\mu_{Leb})$$
with real-valued multiplier function $\lambda \in L^\infty([0,\infty),\mu_{Leb})$.
\end{pro}
\begin{proof}
First note, that $M_\lambda$ is self-adjoint for any real-valued $\lambda$. Furthermore, it is bounded for $\lambda \in L^\infty([0,\infty),\mu_{Leb})$.\newline
Since $M_\lambda$ is injective and non-negative, the function $\lambda$ must be non-zero and non-negative itself. Hence, there exists a $c > 0$ such that the set $S := \{\lambda \geq c\}$ satisfies $\mu_{Leb}\left(S\right)> 0$. On the Hilbert space of all $L^2$-functions supported in $S$, i.e.
\[
\mathcal V := \left\{f \in L^2([0,\infty),\mu_{Leb}) \mid f = f \cdot \chi_S\right\},
\]
the operator $M_\lambda$ is surjective (since for $h \in \mathcal V$ it holds $M_\lambda\left(h \chi_S \lambda^{-1}\right)= h \chi_S = h$) and bounded (since $\chi_S\lambda^{-1} \leq c^{-1} < \infty$). So if $M_\lambda :L^2([0,\infty),\mu_{Leb}) \to L^2([0,\infty),\mu_{Leb})$ was compact, so was $$\text{id}_{\mathcal V} = \left(\left(M_\lambda\right)_{|_\mathcal V}\right)^{_1}\left(M_\lambda\right)_{|_\mathcal V} : \mathcal V \to \mathcal V.$$
However, this is impossible, since $\mathcal V$ is infinite dimensional: There exists a partition $S = \bigcup_{i\in \mathbb N} S_i$ with pairwise disjoint sets $S_i$, $\mu_{Leb}(S_i) > 0$ for all $i \in \mathbb N$, the indicator functions of which are linearly independent.
\end{proof}
}

{\begin{rem}
{\rm Note, that we exploited a specific property of the Lebesgue-measure here, namely that for each set $S \subseteq \R$ with $\mu_{Leb}(S) > 0$ there exists an infinite partition $S = \bigcup_{i\in \mathbb N} S_i$ of pairwise disjoint sets $S_i$ with $\mu_{Leb}(S_i) > 0$ for all $i \in \mathbb N$. However, this property is also true for \emph{any} measure \emph{without atoms}, including e.g. the Cantor measure on $[0,1]$.}
\end{rem}}

\section{Examples} \label{sec:examples}

\subsection{The Hausdorff moment operator} \label{sec:hausop}

\begin{ex}[The Hausdorff moment operator] \label{ex:Haus}
{\rm Apart from the simple multiplication operators on $\R$ outlined in Example~\ref{ex:Frank}, the operator $A: L^2(0,1) \to \ell^2$ defined as
\begin{equation} \label{eq:Haus1}
[A\,x]_j:=\int  _0^1 t^{j-1}\,x(t)\,dt \qquad (j=1,2,...)
\end{equation}
associated with the Hausdorff moment problem (cf.~\cite{Haus23}) is one of the most well-known non-compact bounded linear operators with non-closed range. We consider as operator $T:\X \to \Y$ from \eqref{eq:model} with $\X:=\ell^2$ and $\Y:=L^2(0,1)$ the adjoint operator $T:=A^*:\ell^2 \to L^2(0,1)$ to $A$ from \eqref{eq:Haus1}
of the form
\begin{equation} \label{eq:Haus2}
[T\,u](t):=\sum \limits _{j=1}^\infty u_j\, t^{j-1} \qquad (0 \le t \le 1),
\end{equation}
where the corresponding operator equation $Tu=g$ is also ill-posed of type~I. Then we have, for $T$ from \eqref{eq:Haus2}, a coincidence of $T^*T$ with the infinite Hilbert matrix $\mathcal{H}$ such that $T^*T=(\mathcal{H}_{i,j})_{i,j=1}^{\infty}: \ell^2 \to \ell^2$ takes place with the entries
$$\mathcal{H}_{i,j} = \left(\frac{1}{i + j -1}\right)  \qquad (i,j=1,2,...).$$
It is well-known that $T^*T$ is a non-compact linear bounded operator mapping in $\ell^2$ with  $\|T^*T\|=\pi$ and pure continuous spectrum ${\rm spec}(T^*T)=[0,\pi]$.
For further details concerning the Hausdorff moment problem we also refer, for example, to the studies in \cite{Gerth21}.

Concerning the spectral theorem applied to that non-negative self-adjoint operator $T^*T: \ell^2 \to \ell^2$ we learn from \cite{Rosenblum58} that the measure space $(\Omega,\mu)$ with $\Omega=[0,\infty)$ and the corresponding Lebesgue measure $\mu=\leb$ on $\R$ is appropriate here and yields as multiplier
function $\lambda$ of the multiplication operator $M_\lambda$ from \eqref{eq:Multop} the function
\begin{equation} \label{eq:Haus3}
\lambda(\omega)=\frac{\pi}{\cosh(\pi \omega)} \quad (\omega \in [0,\infty)\;\,\mbox{a.e.}).
\end{equation}
With respect to the function $\lambda$ from \eqref{eq:Haus3} we simply derive as distribution function
\begin{equation} \label{eq:Haus4}
\Phi_{\lambda,\mu_{Leb}}(\epsilon) \approx \frac{1}{\pi}\log\left(\frac{2\pi}{\epsilon} \right)\quad \mbox{for sufficiently small}\;\epsilon>0.
\end{equation}
{Evidently, Assumptions~\ref{ass:main1} and \ref{ass:main2} agree in that example,
and formula~\eqref{eq:Haus4}} gives as in case (b) of Example~\ref{ex:Frank} that $[A_{\Phi_{\lambda,\mu_{Leb}}},B_{\Phi_{\lambda,\mu_{Leb}}}]=\{\infty\}$, which
indicates severe ill-posedness of the Hausdorff moment problem.
}\hfill\fbox{}\end{ex}

\begin{rem}
{\rm If the measure $\mu$ were not fixed to be the Lebesque measure, then by choosing  $\kappa(\omega) =  \frac{1}{2}  \exp\left(\pi \omega\right)$ in Lemma \ref{lem:non_uniqueness}, we obtain a different spectral decomposition with $ \tilde \mu = \kappa\, \leb$, which yields
\begin{align*}
\tilde \Phi_\lambda &= \tilde \mu \left(\{ \omega \in(0,\infty) : \lambda(\omega) > \epsilon\}\right) \\
& \asymp \int_0^{\frac1\pi \log \left(\frac{2\pi}{\epsilon}\right)} \kappa(\omega) \d \omega \\
& = \frac{1}{2\pi}{\exp\left(\pi \omega\right)} \Big|_0^{\frac1\pi \log \left(\frac{2\pi}{\epsilon}\right)}= \frac{1}{\epsilon}-\frac{1}{2\pi} \asymp \frac{1}{\epsilon} \quad \mbox{as} \quad \epsilon \to 0.
\end{align*}
This way, a completely different asymptotic behavior indicating only moderate ill-posedness could be obtained. However, such measure $\tilde \mu$ would in an unnatural way up-weight the high frequencies exponentially.
}\end{rem}

\subsection{Convolution operators} \label{sec:conop}

In this section, let us consider the injective linear convolution operator $T : L^2 \left(\R^d\right) \to L^2 \left(\R^d\right)$, defined by
\begin{equation} \label{eq:convop}
\left[Tu\right](t) = \int_{\R^d} h(t-\tau)\, u(\tau) \,\mathrm d \tau \qquad (t \in \R^d)\,,
\end{equation}
with a convolution kernel $h \in L^1 \left(\R^d\right)$. {The measure space $(\Omega,\mu)$ is here $(\R^d,\mu_{Leb})$ with the multi-dimensional Lebesgue measure over $\R^d\;(d \in \N)$ and $\mu_{Leb}(\Omega)=\infty$. Injectivity of $T$ requires} the condition $$\leb(\{\tau \in \R^d:\,h(\tau)=0\})=0.$$ In this case, it follows from the Fourier convolution theorem that
\[
T = \mathcal F^{-1} M_{\mathcal Fh}\, \mathcal F\,,
\]
with the Fourier transform $[\mathcal F u] (\xi) := \int_{\R^d} \exp\left(- \textup{i} \pi \xi x\right) u(x) \,\mathrm d x\;(\xi \in \R^d)$. Hence a spectral decomposition of $T^*T$ is given as
\[
\mathcal F\, T^*T \mathcal F^{-1} = M_{|\mathcal Fh|^2}\,.
\]
For $h \in L^1 \left(\R^d\right)$, we have that the complex function $[\mathcal F h](\xi)\;(\xi \in \R^d)$ is bounded, continuous and tends to zero as $\|\xi\|$ tends to infinity.
Thus, also the non-negative multiplier function $\lambda(\omega)=|[\mathcal Fh](\omega)|^2 \;(\omega \in \R^d)$ is bounded, continuous and tends to zero as $\|\omega\| \to \infty$. This indicates for $h \not=0$, which is excluded by the injectivity of $T$, pure continuous spectrum of $T^*T$, non-compactness of this operator and,
as a consequence of $0  \in {\rm essran}(\lambda)={\rm spec}(T^*T)$, ill-posedness of type~I in the sense of Nashed of the underlying linear integral equation \eqref{eq:model} of the first kind of convolution type in the Hilbert space $L^2 \left(\R^d\right)$.

\begin{pro}\label{pro:Lpconv}
The linear convolution operator \eqref{eq:convop} mapping in $L^2 \left(\R^d\right)$ is well-defined and bounded as soon as $h \in L^1 \left(\R^d\right)$. Moreover, the multiplier function $\lambda$ assigned to $T^*T$ for a measure space with $\Omega=\R^d$ and the Lebesgue measure $\leb$ over $\R^d$ satisfies Assumptions~\ref{ass:main1} {and~\ref{ass:main2}}
whenever $h \in  L^p \left(\R^d\right)$ for some $1<p \le 2$.
\end{pro}
\begin{proof}
	Let $p \in (1,\infty)$ with $1/p + 1/q = 1$. Then Hölder's inequality implies
\begin{multline*}
	|\left(f * h \right) (x)| \leq \int_{\R^d} \left(|f(y)| |h(x-y)|^{1/q}\right) |h(x-y)|^{1/p}\,\mathrm d y \\\leq \norm{f(\cdot)|h(x-\cdot)|^{1/q}}{L^q} \norm{|h(x-\cdot)|^{1/p}}{L^p} = \norm{f(\cdot)|h(x-\cdot)|^{1/q}}{L^q} \norm{h}{L^1}^{1/p}
\end{multline*}
almost everywhere. Therewith, it follows from Fubini's theorem that
\begin{multline*}
	\norm{f*h}{L^q}^q  \leq \norm{h}{L^1}^{q/p} \int_{\R^d} \norm{f(\cdot)|h(x-\cdot)|^{1/q}}{L^q}^q\,\mathrm d x = \norm{h}{L^1}^{q/p} \int_{\R^d}\int_{\R^d} |f(y)|^q |h(x-y)| \,\mathrm d y \,\mathrm d x \\
	= \norm{h}{L^1}^{q/p} \int_{\R^d}|f(y)|^q \int_{\R^d} |h(x-y)| \,\mathrm d x \,\mathrm d y = \norm{h}{L^1}^{q/p+1} \int_{\R^d}|f(y)|^q \,\mathrm d y = \norm{h}{L^1}^{q}\norm{f}{L^q}^q.
\end{multline*}
This shows, that the convolution operator $T$ in \eqref{eq:convop} is a bounded operator $T : L^q\left(\R^d\right) \to L^q\left(\R^d\right)$ whenever $h \in L^1 \left(\R^d\right)$.

Evidently, we have $\leb(\Omega)=\infty$ for $\Omega=\R^d$ yielding item (a) of Assumption~\ref{ass:main1}. The assertion of item (b) of  Assumption~\ref{ass:main1}, however, is for $1<p \le 2$ an immediate consequence of the Hausdorff-Young inequality
$$\|\mathcal F h\|_{L^q(\R^d)} \le C\,\|h\|_{L^p(\R^d)} \quad \mbox{for} \quad 1 < p \le 2, \quad \frac{1}{p}+\frac{1}{q}=1 \quad \mbox{and} \quad 0<C<\infty $$
(see, e.g, \cite[Theorem~1]{Beckner75}). This implies for  $h \in  L^p \left(\R^d\right)\;(1<p \le 2)$  that we have $\lambda=|\mathcal Fh|^2 \in L^r(\R^d)$ with $r=\frac{p}{2(p-1)} \in [1,\infty)$. Consequently, Corollary~\ref{cor:Lp} is applicable. {The proof is complete, because
Assumption~\ref{ass:main2} is trivially satisfied.}
\end{proof}

\begin{claim}
If a convolution operator $T$ is bounded as $T : L^2 \left(\R^d\right)\to H^s\left(\R^d\right)$ with some $s > 0$, then there exists a $p \in [1,\infty)$ such that the multiplier function $\lambda = |\mathcal F h|^2$ in the spectral decomposition $T^*T = \mathcal F^{-1} M_\lambda \mathcal F$ satisfies $\lambda \in L^p \left(\Omega, \mu\right)$.
\end{claim}
\begin{proof}
Since
\[
\norm{g}{H^s(\R^d)}^2 = \int_{\R^d} \left(1+\|\xi\|^2\right)^s \left| (\mathcal F g)(\xi)\right|^2\d \xi = \norm{ M_\chi \mathcal F g}{L^2 (\R^d)}^2
\]
with $\chi(\xi)  =\left(1+\|\xi\|^2\right)^s, \xi \in\R^d$, the composition
\[
M_\chi \mathcal F  T = M_\chi M_{\mathcal F h} \mathcal F = M_{\chi \cdot \mathcal F h} \mathcal F
\]
is bounded as an operator $L^2 (\R^d) \to L^2 (\R^d)$. Therefore, $M_{\chi \cdot \mathcal F h}:L^2 (\R^d) \to L^2 (\R^d)$ needs to be bounded, which is the case if and only if $\chi \cdot \mathcal F h \in L^\infty\left(\R^d\right)$. This proves the claim.
\end{proof}

\begin{ex}[Simple Gaussian convolution kernel]{\rm Consider the Gaussian convolution kernel $h(t)=\exp(-\|t\|^2)\;(t \in \R^d)$ with the associated multiplier function
$$\lambda(\omega)=|[\mathcal F  h](\omega)|^2=\pi^d\,\exp\left(-\|\omega\|^2/2\right)  \qquad (\omega \in \R^d) .$$
Then we have $\lambda(\omega)> \epsilon$ for $\|\omega\| < R:=\sqrt{2\log(\pi^d/\epsilon)}$. However, the volume of a $d$-ball with radius $R$ is $V_d(R)=\frac{\pi^{d/2}}{\Gamma(d/2)\,d}\,R^d$. This gives
$$ \Phi_{\lambda,\mu_{Leb}}(\epsilon)=\leb\left(\omega \in \R^d:\, \|\omega\| < \sqrt{2\log(\pi^d/\epsilon)} \right)=\frac{\pi^{d/2}}{\Gamma(d/2)\,d}\,\left(\sqrt{2\log(\pi^d/\epsilon)} \right)^d$$
and leads to the asymptotics
$$ \Phi_{\lambda,\mu_{Leb}}(\epsilon)  \asymp \;\left[\log(1/\epsilon)\right]^{\frac{d}{2}} \quad \mbox{as}\quad \epsilon \searrow 0, $$
which shows $A_{\Phi_{\lambda,\mu_{Leb}}}=B_{\Phi_{\lambda,\mu_{Leb}}}=\infty$ indicating severe ill-posedness.
}\end{ex}

	\begin{ex}[Convolution kernel leading to mild ill-posedness]{\rm Consider the convolution kernel $h$ on $\R^d$ defined by means of its Fourier transform
	\[
	\mathcal F h (\xi) = \frac{1}{\left(1+b\norm{\xi}{2}^2\right)^a}
	\]
	for parameters $a,b > 0$. Such Laplace-type kernels have been used in super-resolution microscopy to approximate the point spread function of a STED microscope, see \cite{pwm18}.
	The associated multiplier function $\lambda$ is given by
	$$\lambda(\omega)=|[\mathcal F  h](\omega)|^2=\frac{1}{\left(1+b\norm{\xi}{2}^2\right)^{2a}}  \qquad (\omega \in \R^d),$$
	which implies similarly to Example \ref{ex:Frank} that
			$$ \Phi_{\lambda,\mu_{Leb}}(\epsilon)  \asymp \epsilon^{-\frac1{4a}}\quad \mbox{as}\quad \epsilon \searrow 0, $$
			indicating a degree $2a$ of ill-posedness.
	}\end{ex}

\section{The case of unbounded operators} \label{sec:unbounded}

As in Section~\ref{sec:Halmos} we discuss the application of the spectral theorem, but with one substantial change: the injective linear operator $T: D(T) \subset \X \to \Y$ with domain $D(T)$ dense in $\X$ and hence its self-adjoint non-negative companion $T^*T: D(T) \subset \X \to \X$ with $0 \in {\rm spec}(T^*T)$ are \emph{unbounded} and consequently not compact.
Then the spectral theorem in multiplication operator version is also available, and we refer for details to \cite[Chapt.~VII.3]{w05}. In particular, now in the factorization \eqref{eq:Multop} there occurs an unbounded multiplication operator $M_\lambda$ (cf.~\eqref{eq:lambda}) with a $\mu$-measurable multiplier function $\lambda: \Omega \to \R$ such that $\min {\rm essran}(\lambda)=0$ and $\sup {\rm essran}(\lambda) = +\infty$, since we have again the coincidence ${\rm spec}(T^*T)={\rm essran}(\lambda)$. Let Assumption~\ref{ass:main1} be valid, which implies that $\mu(\Omega)=\infty$, the distribution function $\Phi_{\lambda,\mu_{Leb}}(\epsilon)$ (cf.~\eqref{eq:Phi}) is finite for all $\epsilon>0$ and its inverse $\lambda^*$, the decreasing rearrangement of $\lambda$, is well-defined and an index function at the infinity. It is not difficult to understand
that poles of $\lambda$ do not really influence the decay rate of $\lambda^*(t)$ as $t \to \infty$. Thus, Definition~\ref{def:degree_noncompact} seems to be applicable also here for verifying the degree of ill-posedness. Regularization approaches for ill-posed operator equations with unbounded forward operators
are, for example, discussed in \cite{HMW09} and \cite{Umbricht21} by exploiting the role of corresponding multiplier functions.

\begin{ex}[Fractional integral operators on $\R$] {\rm Following \cite{TauGor99} we consider with $\X=\Y=L^2(\R)$ the family of \emph{fractional integration operators} $T_s: L^2(\R) \to L^2(\R)$. This family depending on the parameter $s>0$ is defined as
\begin{equation} \label{eq:frac}
[T_s(u)](t):=\frac{1}{\Gamma(s)} \int_{-\infty}^t \frac{u(\tau)}{(t-\tau)^{1-s}} \,d\tau \qquad (t \in \R).
\end{equation}
All these operators $T_s$ are densely defined, injective and closed on their domains, which can be introduced by means of the Fourier transform $\mathcal F$ as
$$D(T_s)=\{u \in L^2(\R): |\xi|^{-s}[\mathcal F u](\xi) \in L^2(\R)\}. $$
The corresponding operator $T_s^*T_s$, however, is self-adjoint with domain
$$D(T_s^*T_s)=\{u \in L^2(\R): |\xi|^{-2s}[\mathcal F u](\xi) \in L^2(\R)\}. $$
Moreover, we obtain with respect to the Lebesgue measure on $\R$ the multiplication operator $\mathcal F T^*T \mathcal F^{-1} = M_\lambda$ with the multiplier function
\begin{equation} \label{eq:fraclambda}
\lambda(\omega)=|\omega|^{-2s} \qquad (\omega \in \Omega=\R\; \mbox{a.e.}).
\end{equation}
The inspection of formula \eqref{eq:fraclambda} with some pole at $\omega=0$ shows that, for all $s>0$, ${\rm essran}(\lambda)=[0,\infty)={\rm spec}(T_s^*T_s)$.
Hence the unbounded operator $T_s^*T_s$ has continuous spectrum. In analogy to item (a1) of Example~\ref{ex:Frank} on can show that the problem is moderately ill-posed with degree $s>0$.
}\hfill\fbox{}\end{ex}

\begin{ex}[Source identification in parabolic PDEs] {\rm We consider the final value problem for the non-homogeneous heat equation
$$ \left\{\begin{array}{lcr} v_t(x,t)-\kappa^2\Delta v(x,t)=u(x),& x \in \mathbb{R}^d,& t>0,\\v(x,0)=0,& x \in \mathbb{R}^d, & \\v(x,t_0)=g(x),  &  x \in \mathbb{R}^d,& t_0>0.\end{array}\right. $$
For $\X=\Y=L^2(\R^d)$, given $\kappa>0$ and $t_0>0$ let the equation \eqref{eq:model} represent a source identification problem with the forward operator $T: u \mapsto g$, which is implicitly given by the above final value problem.
Fourier transform $\mathcal F$ with respect to $x$ yields the associated multiplication operator for $T^*T$ of the form $\;\mathcal F\, T^*T \mathcal F^{-1} = M_{\lambda}\;$
with the multiplier function
\begin{equation} \label{eq:PDElambda}
\lambda(\omega)= \frac{\left[1-\exp\left(-t_0\,\kappa^2\|\omega\|_2^2\right)\right]^2}{\kappa^4\|\omega\|_2^4} \qquad (\omega \in \R^d).
\end{equation}
One easily finds that here ${\rm spec}(T^*T)={\rm essran}(\lambda)=[0,\infty)$, which means that $T^*T$ is a non-compact, unbounded linear operator, because $\lambda(\omega)$ has a pole at $\omega=0$. Moreover, ill-posedness takes place with $\lim_{\|\omega\|_2 \to \infty} \lambda(\omega)=0$.
Now we have $\lambda(\omega) \asymp \|\omega\|_2^{-4}$ as $\|\omega\|_2\to \infty$, which should imply for the decreasing rearrangement $\lambda^*$ of $\lambda$  that $\lambda^*(t) \asymp t^{-\frac{4}{d}}$  as $t \to \infty$, which indicates moderate ill-posedness with degree $s=\frac{2}{d}$.
By applying for sufficiently small $\epsilon>0$ the distribution function $\Phi_{\lambda,\mu_{Leb}}(\epsilon) \approx \leb(\{\omega \in \R^d: \|\omega\|_2 \le \frac{1}{\kappa \epsilon^{1/4})}\})$, the denominator $d$ in the asymptotics of $\lambda^*$ arises from the fact that the volume  of
a ball in $\R^d$ with radius $\frac{1}{\kappa \epsilon^{1/4}}$ is proportional to $\epsilon^{-d/4}$.

We refer to \cite{Umbricht21} for more details and the discussion of the full parabolic source identification problem in the multidimensional case.
Let us note here that such degree of ill-posedness $s=\frac{2}{d}$ with dimension $d$ in the denominator also occurs for various compact multivariate forward operators, e.g.~for the inverse of the negative Laplacian in $L^2(\Sigma)$ with homogeneous boundary conditions on a bounded domain $\Sigma \subset \R^d$ and with such $\Sigma$ for the embedding operator from $H^1(\Sigma)$ into $L^2(\Sigma)$ {(see also Example \ref{ex:Laplace})}.
}\hfill\fbox{}\end{ex}

\section{Conclusion and outlook} \label{sec:outlook}

In this paper, we have presented a unified framework for measuring the degree of ill-posedness for linear operator {equations} by means of its spectral decomposition. We have shown that the derived concept coincides with the well-known definition in the compact case and allows us to interpret several non-compact operator {equations} as moderately or severely ill-posed as well {when the Lebesgue measure is chosen as benchmark case. Unfortunately, cross comparisons between compact and non-compact operators cannot be handled by our framework, because compact operators and the Lebesgue measure
are not compatible with respect to multiplication operators, and we refer to Proposition~\ref{pro:nocompact} in this context.}

{Currently our concept is also limited concerning the special situations (A), (B) and (C) in Figure \ref{fig:operators}, which have to be excluded for different reasons.} It will be an interesting topic for future research to investigate whether a situation of finite measure $\mu(\Omega)$ can be handled similarly or not. Furthermore it is unclear if other concepts to measure ill-posedness by means of essential zeros of $\lambda$ in the situation (B) can be related to our concept in some way. Finally, it would be an interesting question to find real-world examples for the situation (C). This mostly seems to be a measure-theoretic problem.

Compared to the compact subcase, it is also completely open what can be said about the degree of ill-posedness of the composition of two operators. The Courant-Fischer theorem yields an upper bound for the decay rate of the singular values of the composition of two compact operators, and hence a lower bound for its degree of ill-posedness. However, in the non-compact case no such result is known yet, and it is not clear if anything can be proven about a composition's degree of ill-posedness with our concept. Furthermore, if a compact and a non-compact operator are composed, strange things can happen as has been recently demonstrated in \cite{DFH24,HofMat22,HofWolf09}.

\section*{Acknowledgment}
The second named author has been supported by the German Science Foundation (DFG) under grant~HO~1454/13-1 (Project No.~453804957).

\section*{Appendix: A glimpse of the finite measure case}

In this appendix, we return to the case of bounded linear operators $T:\X \to \Y$, and we will briefly summarize references and specific approaches for the alternative case of finite measures \linebreak $\mu(\Omega)<\infty$. This case has been comprehensively studied in the literature for non-compact \emph{multiplication operators} $T$ mapping in $L^2(\Omega)$ for some bounded subinterval $\Omega$ of $\R$. The focus of most papers is on the unit interval $\Omega=[0,1]$ and the Lebesgue measure  $\mu=\leb$ over $\R$. For that measure, non-negative multiplier functions $\lambda \in L^\infty(0,1)$ are under consideration that possess essential zeros inside $[0,1]$. This leads to ill-posed situations of the corresponding operator equations \eqref{eq:model}, where the structure of the zeros seem to play a prominent role for the strength of the ill-posedness.  We refer in this context to \cite{Freitag05,Hof06,HofFlei99,HofWolf05,HofWolf09} and recently \cite{mnh22}. In the latter reference, in particular, the total influence of different essential zeros of the multiplier function $\lambda$ on the \emph{increasing rearrangement} of $\lambda$ is discussed.
This increasing rearrangement (see for details \cite{EHZ93}) and its growth rate at a neighbourhood of zero are tools for characterizing the degree of ill-posedness for this class of non-compact operators modelled with finite measure spaces.

Unfortunately, the cross connections of this increasing rearrangement concept for $\Omega=[0,1]$ and the associated decreasing rearrangement concept for $\mu(\Omega)=\infty$, see formula \eqref{eq:invphi} above, are not completely clear.  It is appropriate for the finite measure case to use instead of \eqref{eq:Phi} the non-decreasing distribution function
$$d_\lambda(\epsilon):= \mu \left( \{\omega \in [0,1] :\, \lambda(\omega) \le \epsilon\}\right) \quad (0 \le \epsilon \le \|T\|)$$
of the non-negative multiplier function $\lambda(\omega)\;(0 \le \omega \le 1)$,
and instead of \eqref{eq:invphi} the increasing rearrangement $\lambda^*$ of $\lambda$ defined as
$$\lambda^*(t)=\sup\left\{\epsilon \in [0,\|T\|] :\,\d_\lambda(\epsilon) \le t\right\}\quad (0 \le t \le 1),$$
which is the (generalized) inverse  to the function $d_\lambda$. Both $d_\lambda$ and $\lambda^*$ are index functions at zero, which means that they are non-decreasing functions with positive values for positive arguments and zero limits as the arguments tend to zero.
Then the most authors characterize the ill-posedness by the growth rate of the function $\lambda^*(t)$ in a right neighborhood of $t=0$.  As shown in \cite[Proposition~2]{mnh22} the strongest essential zero of the multiplier $\lambda$ over the interval $[0,1]$
is responsible for the degree of ill-posedness.

\bibliographystyle{plain}
\bibliography{references.bib}

\end{document}